\documentclass[12pt]{amsart}

\usepackage{amssymb,amsmath,amscd,enumerate,verbatim,fullpage}
\usepackage[colorlinks=true,linkcolor=blue,citecolor=blue]{hyperref}
\usepackage{cleveref}
\usepackage[all]{xy}

\usepackage[all]{xy}

\usepackage{color}

\numberwithin{equation}{section}
\theoremstyle{plain}
\newtheorem{theorem}{Theorem}[section]
\newtheorem{lemma}[theorem]{Lemma}
\newtheorem{proposition}[theorem]{Proposition}
\newtheorem{corollary}[theorem]{Corollary}

\newtheorem{question}[theorem]{Question}

\newtheorem*{maintheorem}{Main Theorem}

\theoremstyle{definition}
\newtheorem{definition}[theorem]{Definition}
\newtheorem{quest}[theorem]{Question}

\theoremstyle{remark}

\newtheorem{example}[theorem]{Example}
\newtheorem{rem}[theorem]{Remark}

\newcommand{\depth}{\operatorname{depth}}

\newcommand{\gr}{\operatorname{gr}}
\newcommand{\grade}{\operatorname{grade}}

\newcommand{\Ht}{\operatorname{ht}}

\newcommand{\Min}{\operatorname{Min}}

\newcommand{\Reg}{\operatorname{Reg}}

\newcommand{\Sing}{\operatorname{Sing}}
\newcommand{\Spec}{\operatorname{Spec}}

\newcommand{\ZZ}{{\mathbb Z}}

\def\mm{{\mathfrak m}}
\def\nn{{\mathfrak n}}
\def\pp{{\mathfrak p}}
\def\qq{{\mathfrak q}}

\newcommand{\kk}{\Bbbk}

\newcommand{\what}{\widehat}

\begin{document}
	
	\title[Perturbations of local domains]{On the perturbations of Noetherian local domains}
	
	\author{ Nguyen Hong Duc}
	\address{TIMAS, Thang Long University, Nghiem Xuan Yem, Hanoi, Vietnam}
	\email{duc.nh@thanglong.edu.vn}

	\author{Hop D. Nguyen}
	\address{Institute of Mathematics, Vietnam Academy of Science and Technology, 18 Hoang Quoc Viet, 10072 Hanoi, Vietnam}
	\email{ngdhop@gmail.com; ndhop@math.ac.vn}
	
	\author{Pham Hung Quy}
	\address{Department of Mathematics, FPT University, Hanoi, Vietnam}
	\email{quyph@fe.edu.vn}
	
	\thanks{}

		\dedicatory{Dedicated to Professor Nguyen Tu Cuong on the occasion of his 75th birthday}

\begin{abstract}
We study how the properties of being reduced, an integral domain, and normal, behave under small perturbations of the defining equations of a noetherian local ring. It is not hard to show that the property of being a local integral domain (reduced, normal ring)  is not stable under small perturbations in general. We prove that perturbation stability holds in the following situations: (1) perturbation of being an integral domain for factorial excellent Henselian local rings; (2) perturbation of normality for excellent local complete intersections containing a field of characteristic zero; and (3) perturbation of reducedness for excellent local complete intersections containing a field of characteristic zero, and for factorial Nagata local rings. 
\end{abstract}

	\maketitle

\section{Introduction}
A noetherian local ring $(R,\mm)$ is said to \emph{have isolated singularity} if $R_\pp$ is regular for any non-maximal prime ideal $\pp$. Let $\kk$ be a field and $S = \kk[\![y_1, \ldots, y_s]\!]$ be a power series ring where $s\ge 1$. In 1956, Samuel \cite{S56} showed that for any series $f \in S$  such that $S/(f)$ has isolated singularity, there exists an integer $N\ge 1$ such that for all $\delta \in (y_1,\ldots,y_s)^N,$ there is a local $k$-algebra automorphism of $S$ sending $f$ to $f + \delta$. Thus if $S/(f)$ is a hypersurface ring with isolated singularity, then replacing $f$ by an element in a sufficiently small adic neighborhood of $f$, does not affect the properties of the ring.  Subsequently, many authors have considered the behaviour of various ring-theoretic and module-theoretic properties and numerical invariants under small perturbations; see, e.g. \cite{CS93, DS24, Du22, E74, EH05, HT97, MQS20, QT23, ST96}. We have information about the behaviour under perturbation of Hilbert--Samuel functions \cite{MQS20, QVDT21, ST96}, homology of complexes \cite{E74}, regular sequences \cite{HT97}, positive characteristic singularities \cite{DS24}, among others.
	
Given a noetherian local ring $(R, \mm)$, we say that a ring property $\mathcal{P}$ is {\em stable under perturbations} for $R$ (or {\em perturbation of $\mathcal{P}$ holds} for $R$) if for any regular element $a$ such that $R/(a)$ satisfies $\mathcal{P}$, we also have that $R/(b)$ satisfies $\mathcal{P}$ for all $N\gg 0$ and all $b\in a+\mm^N$, i.e., for all $b$ lying in a sufficiently small $\mm$-adic neighborhood of $a$. For a class $\mathcal C$ of noetherian local rings, we say that $\mathcal{P}$ is stable under perturbations for $\mathcal C$ if it is stable under perturbations for every member of $\mathcal C$. As proved in \cite[Theorem 3.5]{MQS20}, the transformations sending $R/(a)$ to $R/(a+\delta)$ where $\delta$ is a sufficiently small perturbation, preserves the \emph{Hilbert--Samuel function} of $R/(a)$, namely the Hilbert function of the associated graded ring $\gr_\mm(R/(a))$.  Not all ring-theoretic properties behave alike under perturbations and many basic questions remain open. For instance, in a recent work on $F$-singularities in positive characteristics \cite{DS24}, De Stefani and Smirnov showed the stability of $F$-rational singularities, as well as some subclass of $F$-injective singularities under perturbations. On the other hand, they exhibited in \cite[Example 5.3]{DS24} that $F$-purity and (strong) $F$-regularity are not stable under perturbations. It is not clear whether $F$-injectivity always perturbs.
	
In this note, we look at very basic properties of rings: integrality (being a domain), reducedness and normality, and ask how they behave under fine perturbations. As far as we know, this problem was first proposed by L. Duarte in his Ph.D. thesis \cite[Question 5.2.2]{Du22}.  It is not hard to show that these questions generally have a negative answer. In fact, there exists a regular local ring $(R, \mm)$, a non-zerodivisor $a$ such that $R/(a)$ is a domain, however arbitrarily small perturbations of it need not be domains (see \Cref{exam_integralitynotperturb}). Using a classical construction of Nagata, we also produce a normal local domain $R/(a)$ whose arbitrarily small perturbations may fail to be normal (see \Cref{ex_normalitynotperturb} and also \Cref{rem_reducednessnotperturb} on the instability of reducedness). As a consequence, we get negative answers to the last part of \cite[Questions 5.2.1 and 5.2.2]{Du22}. In all these examples, the ring $R$ is regular but not $\mm$-adically complete, so we are naturally led to the following 
\begin{question} 
\label{quest_perturb}
Let $(R, \mm)$ be a complete noetherian local ring, and $a \in \mm$ a regular element. Suppose $R/(a)$ is an integral domain \textup{(}res. reduced, normal ring\textup{)}. Is it true that $R/(a + \delta)$ is an integral domain \textup{(}res. reduced, normal ring\textup{)} for all $N\gg 0$ and all $\delta \in \mm^N$?
\end{question}
We do not know the full answer to this question. On the other hand, using strong Artin approximation, a recent generalization of Samuel's work \cite{GP19}, and various Bertini-type theorems, we prove several results supporting a positive answer for it.
\begin{maintheorem} 
The following statements hold true.
\begin{enumerate}[\quad \rm (1)]
\item Perturbation of integrality holds for excellent factorial Henselian local rings, in particular, for complete regular local rings \textup{(}see  \Cref{thm_integrality} and \Cref{cor_integrality_completeRLR} for details\textup{)}.
\item Perturbation of normality holds for excellent local complete intersections that contain a field of characteristic zero \textup{(}\Cref{thm_normal_CIs}\textup{)}.
\item Perturbation of reducedness holds for excellent local complete intersections that contain a field of characteristic zero  and for factorial Nagata local rings\textup{(}Theorems \ref{thm_reduced_CIs} and \ref{reducedness}\textup{)}.
\end{enumerate}
\end{maintheorem}
Note that \Cref{thm_integrality} answers affirmatively a question in Duarte's thesis \cite[Question 5.2.3]{Du22b}. Note also that in our main results, the base ring $R$ is usually ``nice'' but not necessarily complete. At the moment, we do not know whether perturbation of integrality holds for $\mm$-adically complete local complete intersections, even those that contain a field.

Our results on the perturbation stability of normality and reducedness (\Cref{thm_normal_CIs} and \Cref{thm_reduced_CIs}) are consequences of a more general statement.
\begin{theorem}[= \Cref{thm_SerresRegcond_CIs}]
Let $(R,\mm)$ be an excellent local complete intersection that contains a field of characteristic zero. Let $a\in \mm$ be a regular element such that $R/(a)$ satisfies Serre's regularity condition $(R_\nu)$ for some non-negative integer $\nu$. Then there exists $N\ge 1$ such that for every $\delta \in \mm^N$, $R/(a+\delta)$ also satisfies $(R_\nu)$.
\end{theorem}

Since we will encounter types of rings such as (quasi-)excellent, Nagata, and Henselian ones quite often in this note, we recall their definitions in \Cref{sect_background}. In \Cref{sect_perturbations}, we display examples on the instability of (normal) integral domains, thereby answering in the negative \cite[Question 6.1]{DS24}. We also provide the proofs of our main results in this section, hoping to inspire further study on the tantalizing \Cref{quest_perturb}.
	
\subsection*{Acknowledgment.}  We are grateful to the referee for carefully reading the paper and providing valuable suggestions. A part of this work was done while the authors were visiting  Vietnam Institute for Advanced Study in Mathematics (VIASM) and Institute of Mathematics, Vietnam Academy of Science and Technology (VAST) in 2024. The authors would like to thank these institutes for hospitality and support. We are grateful to Alessandro De Stefani, Linquan Ma, Kazuma Shimomoto, and Ilya Smirnov for many useful comments on a previous draft. The second named author (HDN) was supported by the National Foundation for Science and Technology Development (NAFOSTED) under the grant number 101.04-2023.30. He gratefully acknowledges the support from  the Institute of Mathematics, VAST via grant CSCL01.01/25-26, and from VAST via grant CTTH00.01/25-26.
\section{Background}
\label{sect_background}

Let $\mathcal{P}$ be a given property of noetherian local rings.  We say that the property $\mathcal{P}$ is \emph{stable under perturbations} (or \emph{perturbation of $\mathcal{P}$ holds}) for a noetherian local ring $(R, \mm)$ if for any non-zerodivisor $a \in \mm$ such that $R /(a)$ satisfies $\mathcal{P}$, then $R /(a+\delta)$ satisfies $\mathcal{P}$ for all $N\gg 0$ and all $\delta\in \mm^N$. Similarly, we define when $\mathcal{P}$ is stable under perturbations for a class of (resp. for all) noetherian local rings.  Related to the notion of stability under perturbations, we record the well-known fact that regular sequences perturb; see \cite[Corollary 1]{E74} and \cite{HT97}.
\begin{lemma}
\label{lem_perturb_regularelement}
Let $(R, \mm)$ be a noetherian local ring, and $a\in \mm$ an $R$-regular element. Then there exists an $N\ge 1$ such that for all $\delta\in \mm^N$, $a +\delta$ is $R$-regular.
\end{lemma}

We say that {\em deformation of $\mathcal{P}$ holds} if $R$ satisfies $\mathcal{P}$, as soon as there exists an $R$-regular element $a \in \mm$  such that $R /(a)$ satisfies $\mathcal{P}$. Perturbation and deformation are closely related as shown in \cite[Theorem 2.4]{DS24}. In some sense, perturbation of $\mathcal{P}$ can be regarded as an $\mm$-adic analogue of deformation.  

\begin{theorem}[De Stefani and Smirnov {\cite[Theorem 2.4]{DS24}}]
\label{thm_perturbvsdeform}
Let $\mathcal{P}$ be a property of local rings such that:
\begin{enumerate}[\quad \rm (i)]
\item if a local ring $(R,\mm)$ satisfies $\mathcal{P}$, and $T$ is a variable, then $R[T]_{(\mm, T)}$ satisfies $\mathcal{P}$; and,
\item if $(R, \mm) \rightarrow(S, \nn)$ is faithfully flat, and $S$ satisfies $\mathcal{P}$, then so does $R$.
\end{enumerate}
If perturbation of $\mathcal{P}$ holds for all noetherian local rings, so does deformation of $\mathcal{P}$.
\end{theorem}

Recall that for a non-negative integer $\nu$, a noetherian ring $A$ satisfies \emph{Serre's regularity condition $(R_\nu)$} if $A_\pp$ is a regular local ring for any $\pp\in \Spec(A)$ with $\Ht \pp\le \nu$. The ring $A$ satisfies \emph{Serre's condition $(S_\nu)$} if $\depth A_\pp \ge \min\{\nu, \dim A_\pp\}$ for any prime ideal $\pp\in \Spec(A)$. Recall the following fundamental facts.
\begin{lemma}[See {\cite[§23]{Mat86}} and {\cite[Theorem 11.5]{Eis95}}]
\label{lem_Serrecrit}
Let $A$ be a noetherian ring. The following statements hold.
\begin{enumerate}[\quad \rm (i)]
 \item $A$ is reduced if and only if it satisfies both $(R_0)$ and $(S_1)$.
 \item \textup{(Serre's normality criterion)} $A$ is normal if and only if it satisfies both $(R_1)$ and $(S_2)$.
\end{enumerate}
\end{lemma}

We will employ the next deformation result quite often. For a proof of part (ii), see \cite[Theorem, pp. 439-440]{He}. Part (iii) is contained in \cite[Lemma 0]{BR82}. In view of Lemma \ref{lem_Serrecrit}, both parts (iii) and (iv) are
a direct consequence of part (ii). See also the much more general and sweeping result of Li \cite[Theorem 1.4]{Li24}. 
\begin{proposition}
\label{prop_lifting}
Let $(R,\mm)$ be a local ring and let $a\in \mm$ be a non-zerodivisor. 
\begin{enumerate}[\quad \rm (i)]
\item If $R/(a)$ is an integral domain, then so is $R$.
\item More generally, if $R/(a)$ satisfies Serre's conditions $(R_\nu)$ and $(S_{\nu+1})$, then $R$ also satisfies both $(R_\nu)$ and $(S_{\nu+1})$.
\item If $R/(a)$ is normal, then so is $R$.
\item If $R/(a)$ is reduced, then so is $R$.
\end{enumerate}
\end{proposition}
\begin{rem}
(1) Given a noetherian local ring $(R,\mm)$, we say that a property $\mathcal{Q}$ of principal ideals of $R$ is \emph{stable under perturbations} if any principal $R$-ideal $(a)$ generated by an $R$-regular element,  with the property $\mathcal{Q}$, there exists $N\ge 1$ such that for all $\delta \in \mm^N$, the ideal $(a+\delta)$ also satisfies $\mathcal{Q}$.

In this language, the part concerning integral domains in \Cref{quest_perturb} can be rephrased as follows: \emph{Let $(R,\mm)$ be a complete local ring. Then is the property of being a principal prime ideal stable under perturbations?}

The following example shows that the property of being an irreducible (or primary) principal ideal is not stable under perturbations.

Let $R = \kk[\![x, y]\!]$ be the ring of formal power series over a field $\kk$ with $\mathrm{char}\,\kk \neq 2$. We have $(x^2)$  is an irreducible ideal. However for all $\delta = -y^{2n}, n \ge 1$, as $(x^2 - y^{2n}) = (x + y^{n}) (x - y^{n}),$ $(x^2 - y^{2n})$ is not even primary.

(2) Similarly, one can define the ideal-theoretic version of deformation. Given a noetherian local ring $(R,\mm)$, we say that a property $\mathcal{Q}$ of ideals of $R$ is \emph{stable under deformations} if for any  $R$-ideal $I$ and any element $x\in R$ that is $(R/I)$-regular, such that $I+(x)$ has the property $\mathcal{Q}$, then $I$ also has property $\mathcal{Q}$.  

Thus \Cref{prop_lifting} implies that the property of being a radical (prime) ideal is stable under deformations. On the other hand, the following example shows that the property of being an irreducible (or primary) ideal is not stable under deformations.

Let $R = \kk[\![x, y]\!]$ and $I=(xy)$.  Then $a = x +y$ is a non-zerodivisor on $R/I$. Moreover $I+(x+y)=(x+y,x^2)$ is an irreducible ideal of $R$. However, $I$ is not every primary.
\end{rem}

The authors of \cite{DS24} asked if the converse of \Cref{thm_perturbvsdeform}  also holds true.
\begin{question}[{\cite[Question 6.1]{DS24}}]
\label{ques1}
Let $\mathcal{P}$ be a property of local rings that satisfies the assumptions of Theorem \ref{thm_perturbvsdeform}.  Are deformation and perturbation of $\mathcal{P}$ equivalent?
\end{question}
The following example provides a negative answer to this question.  The example also shows that the integrality property is generally unstable under perturbations.

\begin{example}
\label{exam_integralitynotperturb} 
Let $\kk$ be a field of characteristic $0$. Note that the element $1 + x$ is not a square in $\kk [x]_{(x)}$ since otherwise $1 + x = \dfrac{P^2}{Q^2}$ for some $P,Q \in \kk[x]$, so $2 \deg P=1+2\deg Q,$ which is impossible. However in $\kk [\![x]\!]$ there exists an element $u= 1 + a_1x + a_2 x^2 + \cdots $ such that $u^2 = 1 + x$, thanks to the Taylor expansion of $(1+x)^{1/2}$. Let $R=\kk[x,y]_{(x,y)}$. We have $a = x^2 + x^3 - y^2$ is an irreducible polynomial in $\kk[x,y]$, so $R/(a)$ is a domain. Taking the $(x, y)$-adic completion, we have 
   $$x^2 + x^3 - 
 y^2 = (ux)^2 - y^2 = (ux + y)(ux - y),$$
which is reducible as a power series. For each $N \ge 2$, set
 $$b_{1, N} = x + a_1x^2 + \cdots + a_{N-1}x^N + y \in \kk[x, y],$$
  $$b_{2, N} = x + a_1x^2 + \cdots + a_{N-1}x^N - y \in \kk[x, y].$$
 It is easy to check that $a = b_{1, N} b_{2, N} - \delta_N$ with $\delta_N \in (x, y)^{N+1}R$. Hence for every $N \ge 1$ we have an element $\delta_N \in (x, y)^{N+1}R$ such that $a + \delta_N$ is reducible.
\end{example}
\begin{corollary}
\label{cor_integralitynotperturb}
There exists a property $\mathcal{P}$ of local rings that satisfies the assumptions of Theorem \ref{thm_perturbvsdeform} such that deformation and perturbation of $\mathcal{P}$ are not equivalent.
\end{corollary}
\begin{proof}
Let $\mathcal{P}$ be the property of ``integrality". It is easy to see that $\mathcal{P}$ satisfies the assumptions of Theorem \ref{thm_perturbvsdeform}. It follows from  \Cref{prop_lifting} that deformation of $\mathcal{P}$ holds for local rings, but perturbation of $\mathcal{P}$ does not hold in general for local rings according to \Cref{exam_integralitynotperturb}. 
\end{proof}
We will see from \Cref{ex_normalitynotperturb} that the property of ``normality'' also satisfies the assumptions of Theorem \ref{thm_perturbvsdeform}, but deformation holds for normality while perturbation does not.

Next, let us recall the notions of quasi-excellent, excellent and Nagata rings, following the Stacks Project \cite{StackCA, StackMoreAlg}.
\begin{definition}
Let $\kk$ be a field, and $A$ be a (possibly non-Noetherian) $\kk$-algebra.
\begin{enumerate}
 \item We say that $A$ is \emph{geometrically normal} over $\kk$ if, for every field extension $\kk'/\kk$, the ring $\kk'\otimes_\kk A$ is normal. Equivalently, $A$ is geometrically normal over $\kk$ if $\kk'\otimes_\kk A$ is a normal ring for every finitely generated field extension $\kk'$ of $\kk$.
 \item We say that $A$ is \emph{geometrically regular} over $\kk$ if $A$ is Noetherian, and for every finitely generated field extension $\kk'$ of $\kk$, $\kk'\otimes_\kk A$ is a regular ring.
\end{enumerate}
\end{definition}

If $(R,\mm)$ is a local ring, let $\what{R}$ be the $\mm$-adic completion of $R$, and the \emph{formal fiber} of $R$ over $\pp \in \Spec(R)$ is defined as $\what{R}\otimes_R \kappa(\pp)$, where $\kappa(\pp):=R_\pp/\pp R_\pp$ denotes the residue field of $R_\pp$.
\begin{definition}[G-rings]
A ring $R$ is called a \emph{G-ring} if $R$ is noetherian and for every prime ideals $\qq \subseteq \pp$, the formal fiber of $R_\pp$ over $\qq R_\pp$, namely $\what{R_\pp}\otimes_{R_\pp} \kappa(\qq)$, is geometrically regular over $\kappa(\qq)$. Here $\widehat{R_{\mathfrak p}}$ denotes the $\mathfrak pR_{\mathfrak p}$-adic completion of the local ring $R_{\mathfrak p}$. 
\end{definition}

For a noetherian ring $S$, the \emph{regular locus} of $\Spec(S)$ is 
$$
\Reg(\Spec(S))=\{\qq \in \Spec(S) \mid S_\qq \, \text{is a regular local ring}\}.
$$ 
The \emph{singular locus} $\Sing(\Spec(S))$ of $\Spec(S)$ is the complement of $\Reg(\Spec(S))$ inside $\Spec(S)$.
\begin{definition}
Let $R$ be a noetherian ring, and $X=\Spec(R)$. We say that $R$ is a \emph{J-2} ring if for any finite type $R$-algebra $S$, the regular locus $\Reg(\Spec(S))$ is an open subset of $\Spec(S)$.
\end{definition}
Now we come to the following crucial notions.
\begin{definition}[Quasi-excellent and excellent rings]
Let $R$ be a ring. We say that $R$ is \emph{quasi-excellent} if $R$ is noetherian, a G-ring, and J-2 ring. We say that $R$ an \emph{excellent} ring, if $R$ is quasi-excellent, and universally catenary.
\end{definition}
\begin{rem}
\label{rem_excellentrings}
By \cite[Lemma 15.52.2, Tag 07QU and Lemma 15.52.3, Tag 07QW]{StackMoreAlg}, the following types of rings are excellent:
\begin{enumerate}
 \item fields, $\ZZ$,
\item complete noetherian local rings,
 \item Dedekind domains with fraction field of characteristics zero,
 \item localizations of a finite type algebra over an excellent ring.
\end{enumerate}
\end{rem}
The following lemma follows from the definition of quasi-excellent rings.
\begin{lemma}
\label{lem_SerreRegcondition_completion}
Let $(R,\mm)$ be a quasi-excellent noetherian local ring. Let $\nu\ge 0$ be an integer. Then $R$ satisfies $(R_\nu)$ if and only if $\what{R}$ satisfies $(R_\nu)$.
\end{lemma}
\begin{proof}
Since $R$ is quasi-excellent, it is a G-ring, so the formal fibers $\what{R}\otimes_R \kappa(\pp)$ are regular for all $\pp \in \Spec(R)$. In particular, all the formal fibers satisfy $(R_\nu)$. Since $R\to \what{R}$ is a flat morphism, by the base change result for the condition $(R_\nu)$ \cite[Proposition 2.2.21]{BH98}, $R$ satisfies $(R_\nu)$ if and only if $\what{R}$ does.
\end{proof}
While perturbation of integral domains fails in general,  fortunately, we can prove that perturbation of integrality holds for excellent factorial Henselian local rings (e.g.\,formal power series rings over fields). The key technique is a version of Artin approximation theorem (\cite{A69, HR, PP75,Pos86}). 
\begin{definition}[Henselian rings]
A noetherian local ring $(R,\mm,\kk)$ is called \emph{Henselian} if for every monic polynomial $f\in R[T]$ and every root $a_0\in \kk$ of $\overline{f}:=f+\mm R[T] \in \kk[T]$ such that $a_0$ is not a root of the derivative $\overline{f}'$ (namely $\overline{f}'(a_0)\neq 0$), there exists an $a\in R$ such that $f(a)=0$ and $a_0=a+\mm$.  
\end{definition}
\begin{rem}
Note that by \cite[Lemma 10.153.9, Tag 04GM and Lemma 10.153.10, Tag 06RS]{StackCA}, the class of Henselian local rings contains complete local rings and the rings of algebraic (convergent) power series over (valued) fields.  
\end{rem}

The version of Artin approximation that we will use in this note is
\begin{theorem}[Strong Artin approximation, {\cite[Corollary 3.17]{HR}}]
\label{thm1}
Let $(R,\mm)$ be an excellent Henselian local ring. Assume that $r, s\ge 1$ are integers, $f_1,\ldots,f_r \in R[y]:=R[y_1,\ldots,y_s]$ are polynomials in the indeterminates $y_1,\ldots,y_s$ with coefficients in $R$. Denote $f(y)=(f_1(y),\ldots,f_r(y)) \in R[y]^{\oplus r}$ \textup{(}finite Cartesian product\textup{)}. Then there exists a function $\beta: \mathbb{N} \longrightarrow \mathbb{N}$ with the property:
for every $c \in \mathbb{N}$ and every $\bar{y} \in R^s$ such that $f(\bar{y}) \in (\mm^{\beta(c)})^{\oplus r}$, there exists $ \widetilde{y} \in R^s$ such that $f(\widetilde{y})=\underbrace{(0,0,\ldots,0)}_{\text{$r$ times}}$ and $\widetilde{y}-\bar{y} \in (\mm^c)^{\oplus s}$.
\end{theorem}

\section{Perturbations of (normal) domains and reduced rings}
\label{sect_perturbations}

Our first main result provides a partial positive answer on the perturbation of integrality when $R$ is an excellent factorial Henselian local ring. 
\begin{theorem}
\label{thm_integrality}
Let $(R,\mm)$ be an excellent factorial Henselian local ring. Let $a \in \mm$ be a regular element such that $R/(a)$ is an integral domain. Then there exists $N\ge 1$ such that $R/(a+\delta)$ is a domain for all $\delta\in \mm^N$.
\end{theorem}
\begin{proof}
Since $R$ is a factorial domain, it suffices to show that if $a\in \mm, a\neq 0$ is irreducible then so is every sufficiently small perturbation $a+\delta$ of it.

Since $R$ is excellent Henselian, applying Theorem \ref{thm1} for $f(y_1,y_2)=y_1y_2-a$ in $R[y_1,y_2]$,  there exists a function $\beta: \mathbb{N} \longrightarrow \mathbb{N}$ such that: For all $ c \in \mathbb{N}$, for all $\bar{y}=(\bar y_1,\bar y_2) \in R^2$ such that $f(\bar{y}) \in \mm^{N}$, where $N = \beta(c),$ there exists $\widetilde{y} \in R^2$ such that $f(\widetilde{y})=0$ and $\widetilde{y}-\bar{y} \in (\mm^c)^{\oplus 2}$.

We will show that $a+\delta$ is irreducible for all $\delta\in \mm^N$, where $N:=\beta(1)$. Assume by contradiction that, $a+\delta$ is not irreducible for some $\delta \in \mm^N$. Then
$$ 
a+\delta =\bar y_1 \bar y_2 \quad \text{for some $\bar y_1, \bar y_2\in  \mm$}.
$$
This implies that 
$$
f(\bar y_1, \bar y_2)=\bar y_1 \bar y_2-a =\delta \in \mm^{N}.
$$
Then there exist $\tilde{y}=(\tilde{y}_1,\tilde{y}_2) \in R^2$ such that $f(\tilde{y}_1,\tilde{y}_2)=0$ and $\tilde{y}-\bar{y} \in \mm^{\oplus 2}$. This means that $a$ is reducible, which is a contradiction. This completes the proof.
\end{proof}
An immediate corollary is the following positive answer to a question of Duarte \cite[Question 5.2.3]{Du22b}.
\begin{corollary}
\label{cor_integrality_completeRLR}
Let $R$ be a complete regular local ring, e.g. a power series ring $\kk[\![y_1,\ldots,y_s]\!]$ over a field $\kk$. Let $a \in \mm$ be a regular element such that $R/(a)$ is an integral domain. Then we can find an $N\ge 1$ such that $R/(a+\delta)$ is a domain for every perturbation $\delta\in \mm^N$.
\end{corollary}
\begin{proof}
This follows from \Cref{thm_integrality} since any regular local ring is factorial, and any complete noetherian local ring is Henselian.
\end{proof}
Let us recall the notion of abstract local complete intersections \cite[Definition 2.3.1]{BH98}.
\begin{definition}
\label{defn_completeintersect}
Let $(R,\mm)$ be a noetherian local ring. We say $R$ is a {\em complete intersection} if there exists a regular local ring $(S,\mm_S)$ and an $S$-regular sequence $x_1, \ldots, x_c$ in  $\mm_S$ such that the $\mm$-adic completion of $R$ satisfies $\what{R}\cong S/(x_1,\ldots,x_c)$.
\end{definition}
In view of \Cref{cor_integrality_completeRLR}, it seems to be natural to ask whether integrality perturbs for  $\mm$-adically complete, local complete intersections.
\begin{quest}
Let $(R,\mm)$ be an $\mm$-adically complete, local complete intersection. Let $a\in \mm$ be a regular element such that $R/(a)$ is an integral domain. Is it true that there exists $N\ge 1$ such that for all $\delta \in \mm^N$, $R/(a+\delta)$ is also an integral domain?
\end{quest}
We do not know the answer to this question. On the other hand, we will prove in \Cref{thm_normal_CIs} that normality perturbs for excellent local complete intersections containing a field of characteristic zero.

A natural question is whether perturbation of normality always holds for local rings. We provide a negative answer to this question in the following example, based on a classical construction of Nagata. We are largely influenced by the exposition of Heinzer, Rotthaus and S. Wiegand \cite[Example 4.15]{HRW21}.
\begin{example}[Nagata's example]
\label{ex_normalitynotperturb}
Let $\kk$ be a field of characteristic zero. We start with the power series ring $\kk[\![x,y]\!]$ and its fraction field $Q(\kk[\![x,y]\!])$. Let $\tau\in x\kk[\![x]\!]$
be the following element that is transcendental over $\kk(x,y)$:
\[
\tau=e^x-1=\sum_{n=1}^\infty \frac{x^n}{n!}.
\]
Let $f:=(y+\tau)^2$, and consider the rings
 \[
A=\kk(x,y,f) \cap \kk[\![x,y]\!] \quad \text{(the intersection is taken inside $Q(\kk[\![x,y]\!])$)},
 \]
and $R=A[z]_{(x,y,z)}$. Denote $a=z^2-f\in R$. We claim that the following statements hold true.
\begin{enumerate}
 \item $R$ is a 3-dimensional regular local ring with the unique maximal ideal $\mm=(x,y,z)R$. Moreover, $R$ is not excellent.
 \item $a$ is $R$-regular and $R/(a)$ is a normal domain.
 \item For each $n\ge 1$, there exists $\delta_n \in \mm^n$ such that $R/(a+\delta_n)$ is not a domain, hence not normal.
\end{enumerate}
\begin{proof}
We first note that $A$ is identical to the following ring in \cite[Proposition 6.19]{HRW21}:  
\[
A'=\kk(x,y,f)\cap \kk[y]_{(y)}[\![x]\!].
\]
Indeed, clearly $A'$ is contained in $A$. Conversely, any element  $p\in A$ can be written as a power series in $x$:
\[
p= a_0(y)+a_1(y)x+a_2(y)x^2+\cdots, \quad \text{where $a_i(y)\in \kk[\![y]\!]$ for each $i$}.
\]
To prove $A \subseteq A'$, it remains to establish the next

\textbf{Observation:} We have $a_i(y)\in \kk[y]_{(y)}$ for each $i\ge 0$.

\emph{Proof of the observation}: Indeed, it is harmless to assume that $a_0(y)\neq 0$, as otherwise we may replace $p$ by $p/x \in A$. Note that each element of $k[x,y,f] \subseteq \kk[y][\![x]\!]$ has the form $g(x)=g_0(y)+g_1(y)x+g_2(y)x^2+\cdots$, where $g_i(y)\in \kk[y]$.

Now $p$ belongs to $\kk(x,y,f)$ means that $p=\dfrac{h(x)}{g(x)}$, where $g(x),h(x) \in \kk[y][\![x]\!]$, and $g(x)\neq 0$. Write $h(x)=h_0(y)+h_1(y)x+h_2(y)x^2+\cdots$, where $h_i(y)\in \kk[y]$. Now since $pg(x)=h(x) \in \kk[y][\![x]\!]$, clearing common powers of $x$, we may assume that $g_0(y)\neq 0$. Comparing the two sides of the last equation and arguing by induction, this yields $a_i(y) \in \kk[y]\left[\dfrac{1}{g_0(y)}\right]$ for every $i\ge 0$. In particular, for all such $i$,
\[
a_i(y) \in \kk[y]\left[\dfrac{1}{g_0(y)}\right] \cap \kk[\![y]\!] \subseteq \kk[y]_{(y)},
\]
as desired. Return now to the claims (1)--(3).

For (1): By \cite[Proposition 6.19]{HRW21}, $A=A'$ is a 2-dimensional regular local ring with the maximal ideal $(x,y)A$. Hence $R$ is a 3-dimensional regular local ring with the maximal ideal $(x,y,z)R$. By \cite[Remark 4.16]{HRW21}, $R/(z)\cong A$ is not an excellent ring, so by \Cref{rem_excellentrings}, neither is $R$.

For (2): By \cite[Example 7, pp.~209--211]{Nag62}, 
\[
D=\frac{A[z]}{(z^2-f)A[z]}
\]
is a 2-dimensional normal local domain. In particular, the unique maximal ideal of $D$ is $(x,y,z)D$, hence
\[
R/(a) =\frac{A[z]_{(x,y,z)}}{(z^2-f)A[z]_{(x,y,z)}} \cong D_{(x,y,z)D} \cong D
\]
is a 2-dimensional normal domain. In particular, $a$ is $R$-regular.

For (3): For each $n\ge 1$, write $f=(y_n+x^{n+1}u_n)^2$, where $y_n\in \kk[x,y]\subseteq A, u_n\in \kk[\![x,y]\!]$ are given by
\[
y_n=y+x+\frac{x^2}{2!}+\cdots+\frac{x^n}{n!}, \quad u_n=\sum_{i=0}^\infty\frac{x^i}{(i+n+1)!}.
\]
We note that $f=y_n^2+x^{n+1}u_n(2y_n+x^{n+1}u_n)=y_n^2+x^{n+1}v_n$, where $v_n:=u_n(2y_n+x^{n+1}u_n) \in \kk[\![x,y]\!]$. Since
\[
v_n = \frac{f-y_n^2}{x^{n+1}} \in \kk(x,y,f)\cap \kk[\![x,y]\!]=A,
\]
the element $\delta_n=x^{n+1}v_n$ belongs to $\mm^{n+1}$. Thus
\[
a+\delta_n=z^2-f+x^{n+1}v_n=z^2-y_n^2=(z-y_n)(z+y_n)
\]
is a product of two elements in $\mm$. In particular, $R/(a+\delta_n)$ is not a domain. Being a noetherian local ring, the last ring therefore cannot be normal.
\end{proof}
\end{example}
\begin{rem}
\label{rem_reducednessnotperturb}
Keep the notations as in \Cref{ex_normalitynotperturb}. By \cite[Remark 4.16]{HRW21}, $(A,(x,y)A)$ is a 2-dimension regular local domain with $\what{A}=\kk[\![x,y]\!]$. Moreover, $A/(f)$ is an integral domain. We show that there exist arbitrarily small perturbations of $A/(f)$ that are not reduced.

As in \Cref{ex_normalitynotperturb}, $v_n \in A$, so $\delta_n=x^{n+1}v_n \in (x,y)^{n+1}$. Now 
\[
f-\delta_n=y_n^2 \in (x,y)^2
\]
so $A/(f-\delta_n)$ is not reduced. Thus reducedness is generally unstable under perturbations.
\end{rem}

In view of \Cref{ex_normalitynotperturb} in which $R$ is regular local but not excellent, the most optimistic question about perturbation of normality seems to be
\begin{question}
\label{quest_normlityperturbs}
Let $(R, \mm)$ be an excellent local ring. Suppose $a \in \mm$ is a regular element such that $R/(a)$ is a normal ring. Does there exist an $N\ge 1$ such that for all $\delta \in \mm^N$, the ring $R/(a + \delta)$ is also normal? 
\end{question}
We answer \Cref{quest_normlityperturbs} positively when $R$ is an excellent complete intersection that contains a field of characteristic zero. More generally, we have
\begin{theorem} 
\label{thm_normal_CIs}
Let $(R,\mm)$ be an excellent local complete intersection that contains a field of characteristic zero. Let $a\in \mm$ be a regular element such that $R/(a)$ is a normal domain. Then there exists $N\ge 1$ such that for every $\delta \in \mm^N$, $R/(a+\delta)$ is again a normal domain.
\end{theorem}
Because of Serre's normality criterion (\Cref{lem_Serrecrit}) and the fact that Cohen--Macaulay local rings always satisfy $(S_\nu)$ for every $\nu$, \Cref{thm_normal_CIs} is an immediate consequence of the following stability result for the regularity condition $(R_\nu)$.
\begin{theorem} 
\label{thm_SerresRegcond_CIs}
Let $(R,\mm)$ be an excellent local complete intersection that contains a field $\kk$ of characteristic zero. Let $a\in \mm$ be a regular element such that $R/(a)$ satisfies Serre's condition $(R_\nu)$ for some non-negative integer $\nu$. Then there exists $N\ge 1$ such that for every $\delta \in \mm^N$, $R/(a+\delta)$ also satisfies $(R_\nu)$.
\end{theorem}

There are two main ingredients in the proof of \Cref{thm_normal_CIs}, which will be presented after \Cref{lem_bertini.domain}.  The first result, due to Greuel--Pham \cite{GP19}, generalizes Samuel's theorem \cite{S56} from hypersurfaces to complete intersections over a field. Before stating the result of Greuel--Pham in the generality that is suitable for our purpose, we have the following remark.
\begin{rem}[On two notions of ``isolated singularities'']
Let $\kk$ be a perfect field. Let $I$ be a proper ideal of $S=\kk[\![y_1,\ldots,y_s]\!]$ that is minimally generated by $f_1,\ldots,f_c$. Let $h$ be the height of $I$. Let $I_h\left(\left[\dfrac{\partial f_i}{\partial x_j}\right]\right)$ be the ideal of $h\times h$ minors of the Jacobian matrix $\left[\dfrac{\partial f_i}{\partial x_j}\right]$ of $I$, and let $j(I)=I+I_h\left(\left[\dfrac{\partial f_i}{\partial x_j}\right]\right)$ be the \emph{Jacobian ideal} of $I$. Assume further that $S/I$ is \emph{equidimensional}, namely for every minimal prime $\pp \in \Min(I)$, we have $\dim(S/\pp)=\dim(S/I)$. 

(1) With the above assumptions, for any prime ideal $P\in \Spec(S)$ containing $I$, the following statements are equivalent:
\begin{enumerate}[\quad \rm (i)]
 \item $(S/I)_P$ is a regular local ring;
 \item $P$ does not contain $j(I)$.
\end{enumerate}
So the last statement can be rephrased as
\[
\Sing(\Spec(S/I))=V(j(I)/I).
\]
For the equivalence of statements (i) and (ii), see \cite[Corollary 3.14(2)(iii)]{Ki24} and \cite[Section 4]{Wa94}.

(2) In particular, consider the local ring $(R,\mm):=(S/I, (y_1,\ldots,y_s)/I)$. Then the following statements are equivalent:
\begin{enumerate}[\quad \rm (i)]
 \item $(R,\mm)$ is an \emph{isolated singularity}, namely $\Sing(\Spec(R))\subseteq \{\mm\}$;
 \item $j(I)$ contains a power of $(y_1,\ldots,y_s)$.
\end{enumerate}
Several authors use condition (ii) as the definition of isolated singularity; this is the case with the notion of isolated complete intersection singularity \cite[p.\,198, lines 13--16]{GP19}. So for complete intersections that are quotients of formal power series rings over a perfect field, our notion of ``isolated singularities'' is compatible with that of \cite{GP19}.

(3) If $\kk$ is not perfect, the two notions of isolated singularities are distinct in general. For example, for a prime $p$, let $\kk=\mathbb{F}_p(t)$ and $I=(x_1^p-tx_2^p)\subseteq S=\kk[x_1,x_2]$. Letting $R=S/I$, then $\Sing(\Spec(R))=\{\mm\}$ but $j(I)=I$ is not $(x_1,x_2)$-primary.
\end{rem}
Here is a version of Greuel--Pham's \cite[Theorem 4.6(2)]{GP19} that we will use; note that we require $\kk$ to be a perfect field.
\begin{theorem}
\label{samuel}
Let $\kk$ be a perfect field. Let $I$ be a proper ideal of $S=\kk[\![y_1,\ldots,y_s]\!]$ that is minimally generated by $f_1,\ldots,f_c$.  If either
\begin{enumerate}[\quad \rm (i)]
\item $\dim S/I=0$, or
\item $\dim (S/I)>0$, and $S/I$ is an isolated complete intersection, 
\end{enumerate}
then there exists $N\geq 1$ such that 
$$S /I\cong S /(f_1+\delta_1,\ldots,f_c+\delta_c)$$ for all elements $\delta_1,\ldots,\delta_c\in (y_1,\ldots,y_s)^N$.
\end{theorem}
\begin{rem}
The statement of \Cref{samuel} is the implication (iii) $\Longrightarrow$ (i) in \cite[Theorem 4.6(2)]{GP19}, which always valid even without $\kk$ being infinite by the proof of  \cite[Theorem 4.6(2)]{GP19} (see, \cite[Problem 4.7]{GP19}).
\end{rem}

The second main ingredient in the proof of \Cref{thm_SerresRegcond_CIs} is the following result on the relationship between a ring $R$ and a relative hypersurface $R[t_1,\ldots,t_n]/(F)$, where $F\in R[t_1,\ldots,t_n]$. In such a situation, we denote by $I_F$ the ideal of $R$ generated the coefficients of $F$. We are mostly interested in the case $F=t_1x_1+\cdots+t_nx_n$, where $x_1,\ldots,x_n \in R$ generate an $\mm$-primary ideal, in which $I_F=(x_1,\ldots,x_n)$. The study of such relative hypersurfaces has been taken up by  Hochster \cite{Ho73}. The following statement is a Bertini-type theorem for the $(R_s)$ property, proved by Trung \cite[Page 222--223]{Tr79} and rediscovered in \cite[Lemma 10]{MS23}.
\begin{lemma}
\label{bertini}
Let $(R, \mm)$ be a noetherian local ring satisfying $(R_\nu)$ such that $\dim R\ge \nu+2$. If $\mm=\left(x_1, \ldots, x_n\right)$, then $\dfrac{R(t_1, \ldots, t_n)}{\left(t_1 x_1+\cdots+t_n x_n\right)}$ still satisfies $(R_\nu)$. Here $R(t_1, \ldots, t_n)$ denotes the localization of the ring $R[t_1, \ldots, t_n]$ at the ideal $\mm R[t_1, \ldots, t_n]$.
\end{lemma}

We also record a criterion for a relative hypersurface section to be reduced and irreducible.

\begin{lemma}[{\cite[Korollar 2.3]{Tr79}}]
\label{lem_bertini.domain.general}
Let $R$ be a noetherian domain. Let $n\ge 1$ an integer, $F\in R[t_1,\ldots,t_n]$ such that $\deg F\ge 1$. Denote by $Q(R)$ the field of fractions of $R$. Then the following statements are equivalent: 
\begin{enumerate}[\quad \rm (i)]
 \item $(F)\subseteq R[t_1,\ldots,t_n]$ is a prime ideal;
 \item $\grade(I_F, R)\ge 2$ and $F$ is irreducible as a polynomial in $Q(R)[t_1,\ldots,t_n]$.
\end{enumerate}
Here, we use the convention that $\grade(R,R)=+\infty$.
\end{lemma}
As a direct consequence of \Cref{lem_bertini.domain.general}, we get
\begin{lemma}
\label{lem_bertini.domain}
Let $(R,\mm)$ be noetherian local domain, where $\mm \neq (0)$. Assume that the elements $x_1,\ldots,x_n \in R$ generate an $\mm$-primary ideal. Then $\dfrac{R[t_1,\ldots,t_n]}{(t_1x_1+\cdots+t_nx_n)}$ is an integral domain if and only if $\depth R\ge 2$.
\end{lemma}

\begin{proof}
Note that for $F=x_1t_1+\cdots+x_nt_n$, $I_F=(x_1,\ldots,x_n)$ is $\mm$-primary so $\grade(I_F,R)=\grade(\mm,R)=\depth R$. Since $\mm \neq (0)$, $x_1,\ldots,x_n$ are not all zero, so $\deg F\ge 1$. The conclusion follows from \Cref{lem_bertini.domain.general} since as a linear form in $Q(R)[t_1,\ldots,t_n]$, $F$ is always irreducible.
\end{proof}

We are ready to present the
\begin{proof}[Proof of Theorem \ref{thm_SerresRegcond_CIs}]
In the first step, we reduce to the case $\dim R \le \nu+2$. In fact, assume that $\dim R\ge \nu+3$. Let $\mm=(x_1,\ldots,x_n)$ and let $F:=t_1x_1+\cdots+t_n x_n$ for new variables $t_1,\ldots,t_n$. Since $R/(a)$ is Cohen--Macaulay of dimension $\dim R-1\ge \nu+2$, $(F)$ is a prime ideal of $(R/(a))(t_1,\ldots,t_n)$ by \Cref{lem_bertini.domain}. In particular, $a,F$ is a regular sequence on $R(t_1,\ldots,t_n)$.  Since $R/(a)$ satisfies $(R_\nu)$, by the local Bertini theorem  \Cref{bertini}, the ring
$$\left(R(t_1,\ldots, t_n)/(F)\right)/(a)\cong \left(R/(a)(t_1,\ldots, t_n)\right)/(F)$$
satisfies $(R_\nu)$ as well.

\textbf{Claim:} $R(t_1,\ldots,t_n)/(F)$ is again a local complete intersection.

Indeed, note that $(R,\mm) \to (R(t_1,\ldots,t_n), \mm R(t_1,\ldots,t_n))$ is flat local morphism with closed fiber $R(t_1,\ldots,t_n)/\mm R(t_1,\ldots,t_n) \cong Q(k[t_1,\ldots,t_n])$, which is a field. Hence as $R$ is a complete intersection, so is $R(t_1,\ldots,t_n)$ by work of Avramov \cite[page 30, (1.9.2)]{Avr77}.

Since $a,F$ is an $R(t_1,\ldots,t_n)$-regular sequence, and the latter is a local ring,  $F$ is regular on $R(t_1,\ldots,t_n)$. Now using \cite[Proposition 3.10]{Avr77} (or \cite[Theorem 2.3.4]{BH98}), $R(t_1,\ldots,t_n)/(F)$ is also a complete intersection, as claimed.

Therefore $R(t_1,\ldots,t_n)/(F)$ is again a local complete intersection containing $\kk$, that is excellent thanks to \Cref{rem_excellentrings}. If we have proved that the ring $\dfrac{R(t_1,\ldots, t_n)/(F)}{(a+\delta)}$ satisfies $(R_\nu)$, for all small $\delta$, then being Cohen--Macaulay, the last ring satisfies both $(R_\nu)$ and $(S_{\nu+1})$. Then \Cref{prop_lifting}(ii) implies that the conditions $(R_\nu)$ and $(S_{\nu+1})$ are also fulfilled by $\left(R/(a+\delta)\right)(t_1,\ldots, t_n)$. Descending along the flat morphism $R/(a+\delta) \to \left(R/(a+\delta)\right)(t_1,\ldots, t_n)$, we deduce that $R/(a+\delta)$ fulfills $(R_\nu)$. Therefore we may replace $R$ by $R(t_1,\ldots, t_n)/(F)$ and reduce the dimension of $R$.

Assume that $\dim R \le \nu+2$ and $R/(a)$ satisfies $(R_\nu)$. Note that being a quotient of an excellent ring, $R/(a)$ is excellent by \Cref{rem_excellentrings}. It follows from \Cref{lem_SerreRegcondition_completion} that $\widehat{R}/(a)$ also satisfies $(R_\nu)$. Since $\what{R}$ contains a field, its coefficient ring is its residue field $K=\what{R}/\mm \what{R} \cong R/\mm$.  By the Cohen structure theorem for local rings containing a field \cite[Theorem 7.7]{Eis95}, there is an isomorphism
$$
\widehat{R} \cong K[\![y_1,\ldots,y_s]\!]/I,
$$
for some ideal $I$ of $S=K[\![y_1,\ldots,y_s]\!]$. Since $R$ is a complete intersection, we may assume that $I$ is generated by an $S$-regular sequence $f_1,\ldots,f_c$ thanks to \cite[Theorem 2.3.3]{BH98}. By abuse of notations, we regard $a\in R\subseteq \what{R}$ as an element of $S$. Then the ring
$$\widehat{R}/(a) \cong K[\![y_1,\ldots,y_s]\!]/(a,f_1,\ldots,f_c)$$ satisfies condition $(R_\nu)$ and has dimension $\dim R-1 \le \nu+1$. This implies that it has at worst isolated singularity. Containing $\kk$ as a subfield, $K\cong R/\mm$ has characteristic zero, and hence is perfect, so applying Theorem \ref{samuel} to the ideal $(a,f_1,\ldots,f_c)$, we get an $N\geq 1$ such that
$$K[\![y_1,\ldots,y_s]\!]/(a,f_1,\ldots,f_c)\cong K[\![y_1,\ldots,y_s]\!]/(a+\delta_0,f_1+\delta_1,\ldots,f_c+\delta_c)$$
for all $\delta_0,\delta_1,\ldots,\delta_c\in (y_1,\ldots,y_s)^N$.
We now take any $\delta\in \mm^N R$, which again can be seen as an element of $S$. Then
$$\widehat{R}/(a+\delta) \cong K[\![y_1,\ldots,y_s]\!]/(a+\delta,f_1,\ldots,f_c)\cong K[\![y_1,\ldots,y_s]\!]/(a,f_1,\ldots,f_c)\cong \widehat{R}/(a),$$
which yields the fulfillment of $(R_\nu)$ of $\widehat{R}/(a+\delta)$. Since $\mm$-adic completion is a flat extension, the ring $R/(a+\delta)$ itself fulfills $(R_\nu)$ for every $\delta\in \mm^N$. The proof is concluded.
\end{proof}
\begin{rem}
	\begin{enumerate}[\quad \rm (i)]
		\item We had hoped to extend the argument to perfect fields of characteristic $p>0$, but we were not able to do so. The difficulty is that, in our proof, we pass from $R$ to $R(t_1,\ldots,t_s)$. Even if the residue field $k$ is perfect, the resulting residue field $k(t_1,\ldots,t_s)$ is not perfect. Thus, the Jacobian criterion in the form used in our argument cannot be applied directly.
		\item The hypothesis that $R$ is excellent in \Cref{thm_normal_CIs} cannot be dropped, as we have seen in \Cref{ex_normalitynotperturb}, where $R$ is a regular local ring that is not excellent. 
		
	\end{enumerate} 
	
\end{rem}

As a consequence of \Cref{lem_Serrecrit} and \Cref{thm_SerresRegcond_CIs}, we get that reducedness perturbs for excellent local complete intersections that contain a field of characteristic zero.
\begin{theorem} 
\label{thm_reduced_CIs}
Let $(R,\mm)$ be an excellent local complete intersection containing a field of characteristic zero. Let $a\in \mm$ be a regular element such that $R/(a)$ is reduced. Then there exists $N\geq 1$ such that for every $\delta \in \mm^N$, $R/(a+\delta)$ is also reduced.
\end{theorem}
\begin{proof}
The ring $R/(a)$ is Cohen--Macaulay, so it always satisfies $(S_\nu)$ for every $\nu$. Thus $R/(a)$ is reduced if and only if it satisfies $(R_0)$. It remains only to apply \Cref{thm_SerresRegcond_CIs} to get the conclusion.
\end{proof}

Our next theorem claims that reducedness is stable under perturbations when $R$ is factorial Nagata. An element $x$ of a local  ring $(R,\mm)$ is \emph{square-free} if it cannot be written as $x=y^2z$ for elements $y, z\in R$ where $y\in \mm$. 
\begin{rem}
\label{rem_squarefree}
Note that if $(R,\mm)$ is a factorial local ring, then for $x\in \mm$, the quotient $R/(x)$ is reduced if and only if $x$ is square-free. 
\end{rem}
Here comes our last main result, which is concerned with stability of reducedness for Nagata rings. 
Recall that an {\em N-2 ring} (or {\em Japanese ring}) is an integral domain $R$ such for every finite field extension $L$ of the fraction field $Q(R)$ of $R$, the integral closure of $R$ in $L$ is $R$-finite.  A ring $R$ is said to be {\em  Nagata} if $R$ is noetherian and for every prime ideal $\pp$ the ring $R/\pp$ is N-2 (\cite[Definition 10.162.1]{StackCA} and \cite[Definition 31.A]{Mat80}). 
\begin{rem}
\label{rem_Nagata}
By \cite[Chapter 12]{Mat80} and \cite[Section 162]{StackCA}, examples of Nagata rings include: 
\begin{enumerate}
\item fields, $\mathbb{Z}$, 
\item complete noetherian local rings,
\item Dedekind domains with fraction field of characteristic zero,
\item localizations of a finite type algebra over a Nagata ring.
\end{enumerate}
\end{rem} 
\begin{theorem}
\label{reducedness}
Reducedness is stable under perturbations for factorial Nagata rings. More precisely, let $(R,\mm)$ be a factorial Nagata ring and let $a\in\mm$ be a nonzero element such that $R/(a)$ is reduced. Then there exists $N\geq1$ such that $R/(a+\delta)$ is reduced for every $\delta\in\mm^N$.
\end{theorem}
\begin{proof}

Let us denote by $\what{R}$ the $\mm$-adic completion of $R$ and let $\what{\mm}:=\mm\what{R}$. Consider the polynomial
$$f(y)=f(y_1,y_2)=a-y_1 y_2^2$$  in $\what{R}[y_1,y_2]$.
By Theorem \ref{thm1}, there exists a function $\beta: \mathbb{N} \longrightarrow \mathbb{N}$ such that: For all $ c \in \mathbb{N}$, for all $\bar{y}\in \what{R}^2$ such that $f(\bar{y}) \in \what{\mm}^{\beta(c)},$ there exists $\widetilde{y} \in \what{R}^2$ such that $f(\widetilde{y})=0$ and $\widetilde{y}-\bar{y} \in (\what{\mm}^c)^{\oplus 2}$.
	
It suffices to show that for all $\delta\in\mm^N$, $N:=\beta(1)$ the ring $R/(a+\delta)$ is reduced, equivalently, $a+\delta$ is square-free in $R$ according to \Cref{rem_squarefree}. Indeed, assume by contradiction that $a+\delta$ is not square-free for some $\delta\in\mm^N$. Then there exist $\bar y_1$ and $\bar y_2$ in $R$  such that $a+\delta=\bar y_1\bar y_2^2$ and $\bar y_2\in \mm$. This implies that $f(\bar y_1, \bar y_2)\in \what{\mm}^{N}.$
Then there exists $\widetilde{y}=(\tilde{y}_1,\tilde{y}_2)\in\what{R}^2$ such that $f(\widetilde{y})=0$ and $\widetilde{y}-\bar{y} \in \what{\mm}^{\oplus 2}$. It yields that $a=\widetilde{y}_1\widetilde{y}_2^2$ and $\widetilde{y}_2\in \what{\mm}$ in $\what{R}$, therefore $\what{R}/(a)$ is not reduced. Since $R$, and hence $R/(a)$, is Nagata ring (\Cref{rem_Nagata}), it follows from \cite[Theorem 70]{Mat80} (or \cite[Lemma 15.44.6, Tag 07NZ]{StackMoreAlg}) that $R/(a)$ is not reduced, a contradiction. This completes the proof.
\end{proof}

\bibliographystyle{amsalpha}

\end{document}